\documentclass{amsart}
\usepackage{amssymb,latexsym,amsmath,amsxtra,graphicx}

\newtheorem{algorithm}{Algorithm}[section]
\newtheorem{theorem}{Theorem}[section]
\newtheorem{lemma}[theorem]{Lemma}

\newcommand{\Z}{\mathbb{Z}}
\newcommand{\R}{\mathbb{R}}

\begin{document}

\title[Inverses and factoring]{The distribution of solutions to $xy = N \mod a$ with an 
application to factoring integers}

\author[Michael Rubinstein]{Michael Rubinstein\\
    \textsl{Pure Mathematics\\University of Waterloo\\200 University Ave W\\Waterloo, Ontario, N2L 3G1\\Canada}
}

\subjclass[2000]{11Y05; 11L050}
\keywords{factoring, number theory, Kloosterman sums}

\begin{abstract} 
We consider the uniform distribution of solutions $(x,y)$ to
$xy=N \mod a$, and obtain a bound on the second moment of the number of
solutions in squares of length approximately $a^{1/2}$. We use this to study a
new factoring algorithm that factors $N=UV$ provably in $O(N^{1/3+\epsilon})$
time, and discuss the potential for improving the runtime to sub-exponential.
\end{abstract}

\maketitle

\section{Introduction}
\label{sec:intro}

Let $\gcd(a,N)=1$. A classic application of Kloosterman sums shows that the
points $(x,y) \mod a$ satisfying $xy = N \mod a$ become uniformly distributed
in the square of side length $a$ as $a \to \infty$. In this paper we
investigate an application of this fact to the problem of factoring integers.
We give a new method to factor the integer $N$ which beats trial division, and
prove that it runs in time $O(N^{1/3+\epsilon})$.

While the complexity of our method is not exciting, considering the existence
of several probabilistic sub-exponential factoring algorithms, the runtime here
is provable and does compete favourably with the best known provable factoring
algorithm, Pollard-Strassen, which only runs in time $O(N^{1/4+\epsilon})$.
Shank's class group method runs in time $O(N^{1/5+\epsilon})$ assuming the GRH.
Our algorithm is described in Section~\ref{sec:algorithm}.

Furthermore, proving this runtime requires understanding the finer distribution of solutions
to $xy = N \mod a$, and our results in this regards are interesting in their own right. 
We discuss the problem on uniform distribution in Sections~\ref{sec:dist} and~\ref{sec:moment}.

Finally, all existing sub-exponential factoring algorithms have grown out of much
weaker exponential algorithms, and we hope that the factoring ideas presented
here will be improved. In Section~\ref{sec:subexp} we discuss some needed improvements
to achieve a better runtime.

We have not implemented the algorithms described in this paper. The purpose
of this paper is to present a new approach to factoring integers and analyse its runtime.

\maketitle

\section{Algorithm- hide and seek}
\label{sec:algorithm}

Let $N$ be a positive integer that we wish to factor.
Say $N=UV$ where $U$ and $V$ are positive integers, not necessarily prime,
with $1 <U \leq V$. For simplicity, assume $V<2U$, so that $V<(2N)^{1/2}$.
The general case, without this restriction, will be handled at the end of this section.

The idea behind the algorithm is to perform trial division
of $N$ by a couple of integers, and to use information about the remainder to
determine the factors $U$ and $V$.

Let $a$ be a positive integer, $1<a<N$. By the division algorithm, write 
\begin{eqnarray}
    \label{eq:linear}
    U=u_1 a + u_0, \quad \text{with $0 \leq u_0 < a$} \notag \\
    V=v_1 a + v_0, \quad \text{with $0 \leq v_0 < a$}. 
\end{eqnarray}
Assume that $u_0$ is relatively prime to $a$, and likewise for $v_0$, since 
otherwise we easily extract a factor of $N$ by taking $\gcd(a,N)$.
If, for a given $a$, we can determine $u_0,u_1,v_0,v_1$ then we have 
found $U$ and $V$. 

Consider $N = u_0 v_0 \mod a$. One cannot simply determine $u_0$
and $v_0$ from the value of $N \mod a$, because $\phi(a)$ pairs of integers
$(x,y) \mod a$  satisfy $xy = N \mod a$
(if $x = m u_0 \mod a$, then $y = m^{-1} v_0 \mod a$, where $\gcd(m,a)=1$).

However, say $a$ is large,  $a \geq \lceil (2N)^{1/3} \rceil > V^{2/3}$, so
that $v_1$ and $u_1$ are comparatively small,
$u_1, v_1\leq V^{1/3}$, i.e. both are $< a^{1/2}$.
If we consider $N \mod a-\delta$
\begin{equation}
    N = UV = (u_1 \delta + u_0)(v_1 \delta + v_0) \mod a-\delta,
\end{equation}
for $\delta=0,1$, we get, as solutions $(x,y)$ to $xy=N \mod a-\delta$, two nearby
points,  $(u_0,v_0)$ and $(u_0+u_1,v_0+v_1)$, whose coordinates are within $a^{1/2}$ of one another.
This pair of points is just one pair amongst the many pairs of solutions to the above equations,
for $\delta=0,1$. However, the fact that the solutions are nearby reduces the amount of checking that we
need to do in order to find the pair of points, $(u_0,v_0)$ and $(u_0+u_1,v_0+v_1)$, that we seek.

Figures~\ref{fig:4 ex} and~\ref{fig:superimposed} illustrate this fact, for $N=
1910861 = 1061 \times 1801$, and $a=157$. Thus, $U=1061$, $u_0=119$, $u_1=8$, and
$V=1801$, $v_0=74$, $v_1=15$. Rather than just depict the solutions to
$xy=N \mod a-\delta$, for $\delta=0,1$, we also plot the solutions for
$\delta=2,3$ (though our algorithm below only makes use of solutions for
$\delta=0,1$). Plotting four sets of solutions, for $\delta=0,1,2,3$ makes it
easier for the human eye to tell the points $(u_0,v_0)=(119,74)$,
$(u_0+u_1,v_0+v_1)=(125,85)$, $(u_0+2u_1,v_0+2v_1)=(131,96)$, and
$(u_0+3u_1,v_0+3v_1)=(137,107)$ from the random coincidences of nearby points
as these all lie equally spaced apart and on one line.

\renewcommand{\thefigure}{\arabic{figure}}

\newpage
\thispagestyle{empty}
\begin{figure}[h]
    \centerline{
       \includegraphics[width=3in,height=3in]{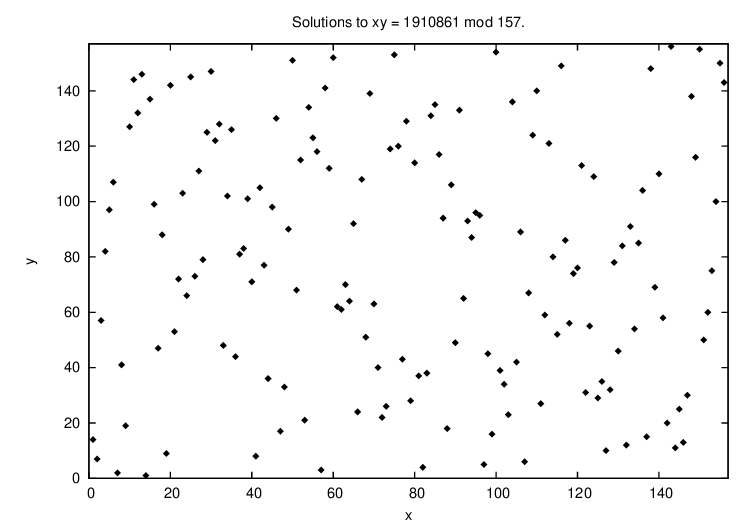}
       \includegraphics[width=3in,height=3in]{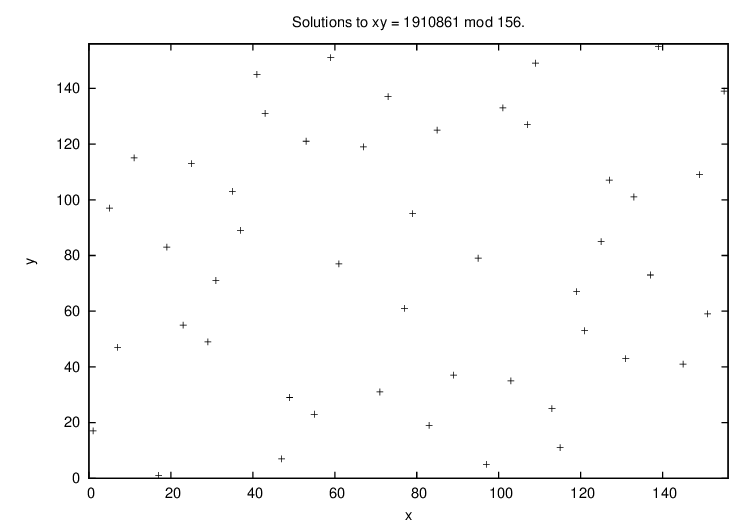}
    }
    \centerline{
       \includegraphics[width=3in,height=3in]{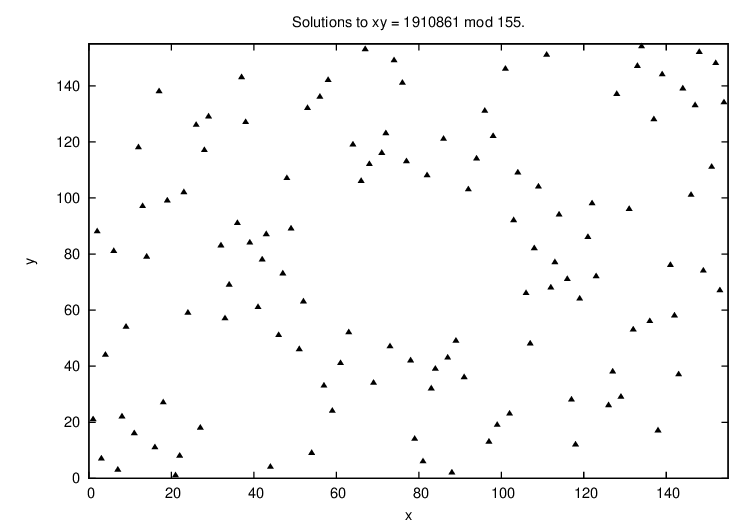}
       \includegraphics[width=3in,height=3in]{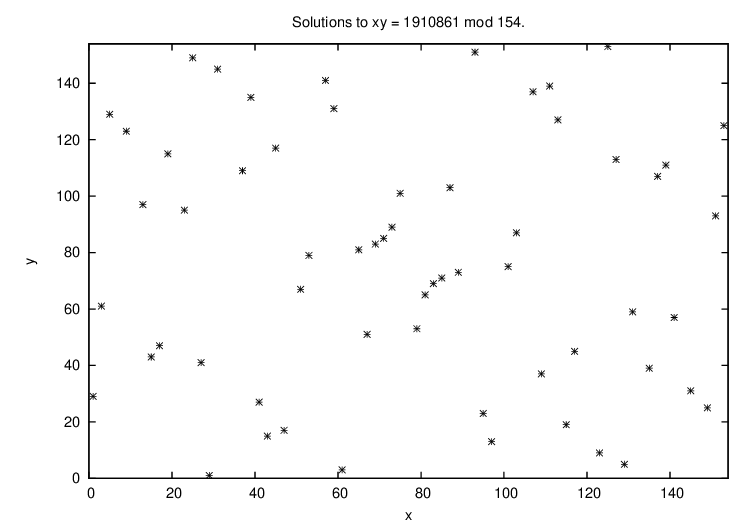}
     }
     \caption[Four solution sets]{The solutions $(x,y)$ to $xy = 1910861 \mod 157-\delta$, for 
     $\delta=0,1,2,3$
     }\label{fig:4 ex}
\end{figure}

\newpage
\thispagestyle{empty}
\begin{figure}[h]
    \centerline{
       \includegraphics[height=5in,width=5in]{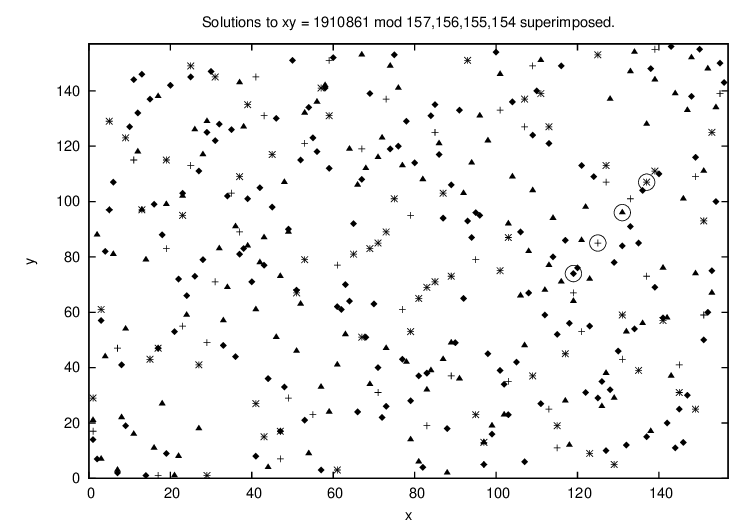}
     }
     \caption[Solution sets superimposed]{The solutions $(x,y)$ to $xy = 1910861 \mod 157-\delta$, for
          $\delta=0,1,2,3$, superimposed. The four circled points are
          $(u_0,v_0)=(119,74)$, $(u_0+u_1,v_0+v_1)=(125,85)$, $(u_0+2u_1,v_0+2v_1)=(131,96)$, and
          $(u_0+3u_1,v_0+3v_1)=(137,107)$.
     }\label{fig:superimposed}
\end{figure}

So, we can set, say, $a=\lceil (2N)^{1/3} \rceil$, and partition the Cartesian
plane into squares of side length $a^{1/2}$, each square being of the form
$\{(x,y)\in\R^2| ma^{1/2} \leq x < (m+1)a^{1/2}, na^{1/2} \leq y < (n+1)a^{1/2} \}$,
where $m,n\in \Z$.

We then list
{\it all} $\phi(a)$ pairs of integers $(x,y)$, with $1 \leq x,y \leq
a$, that satisfy $xy = N \mod a$, throwing them into our
squares of side lengths $a^{1/2}$. We can assume that
$\gcd(a,N)=1$, because, otherwise we easily extract a factor of $N$.

We can compute all inverses mod $a$, and hence all $(x,y)=(x,x^{-1}N) mod a$ in
$O(a)$ operations mod $a$. To compute all inverses, start with $m=2$, multiply
$\mod a$ by $m$ until we arrive at 1, or hit a residue class already
encountered (in which case $m$ is not invertible). Then, take the first residue not yet encountered
and repeat the previous step until all residue classes are exhausted.

Having produced all solutions for the modulus $a$, we then repeat the process
for the modulus $a-1$. For each solution $(x_1,y_1)$ to $xy=N \mod a-1$, we
determine which $a^{1/2} \times a^{1/2}$ square it falls within, and consider
all nearby (with each coordinate within $a^{1/2}$, wrapping to the opposite
side of the larger $a\times a$ square if needed) solutions $(x_0,y_0)$ to $xy =
N \mod a$ from our list of stored solutions. We set $\mu_0=x_0$, $\nu_0=y_0$,
$\mu_1 = x_1-\mu_0$, $\nu_1 = y_1-\nu_0$, and check whether $(\mu_1
a+\mu_0)(\nu_1 a + \nu_0) = N$. If so, we have determined a non-trivial factor
of $N$ and quit.

%

How much work does comparing pairs of points $(x_0,y_0)$ and $(x_1,y_1)$ entail?
There are $\phi(a-1)$ solutions to $xy = N \mod a-1$, and, typically, we expect
there to be only a handful of solutions to $xy=N \mod a$ whose coordinates are
each within $a^{1/2}$. Each such pair of solutions gives us candidate values
$\mu_0,\nu_0$ and $\mu_1,\nu_1$ for $u_0,v_0$ and $u_1,v_1$, and we check to
see whether they produce $N=(u_1 a + u_0)(v_1 a + v_0)$.
On average, each $a^{1/2} \times a^{1/2}$ square
contains $O(1)$ points, the overall time to check all squares and
points is roughly predicted to be $O(a)$. In Section~\ref{sec:dist} we obtain a
runtime bound of $O(a^{1+\epsilon})$. This algorithm terminates successfully when the
true points $(u_0,v_0)$ and $(u_0+u_1,v_0+v_1)$ are found. Since $a=O(N^{1/3})$
this gives a running time that is provably $O(N^{1/3+\epsilon})$.

The idea that lies behind the algorithm suggests the name `Hide and Seek'. The
solutions that we seek $(u_0,v_0)$ and $(u_0+u_1,v_0+v_1)$ are hiding amongst
many solutions in the large $a \times a$ square, but, like children who have
hidden next to one another while playing the game Hide and Seek, they have
become easier to spot.

We summarize the above in the following algorithm.

\begin{algorithm}[Hide and Seek]
    Let $N=UV$ be a positive integer, and assume that $1<U \leq V < 2U$, with $U,V \in \Z$.
    Thus $V<(2N)^{1/2}$.
    For given positive integers $N,r$, define
    \begin{equation}
        H_{N,a} = \left\{ (x,y) \,|\,\, xy = N \mod a, \,\, 0 \leq x,y < a \right\}.
        \label{eq:H}
    \end{equation}
\begin{itemize}
    \item[Step 1] Set $a=\lceil (2N)^{1/3} \rceil$.
    \item[Step 2]Use the Euclidean algorithm to compute $\gcd(N,a-\delta)$ for $\delta=0,1$.
        If either gcd is $>1$ then we have determined a non-trivial factor of $N$ and quit.
    \item[Step 3] Compute and store in an array all $\phi(a)$ points of
        $H_{N,a}$. This can be done using $O(a)$ arithmetic operations mod $a$
        as described above.
    \item[Step 4] For $0 \leq m,n < a^{1/2}$, initialize a doubly indexed
        array, `Bin'. Each element, Bin[m,n], will contain a list of points
        and be used to partition $H_{N,a}$. Each is initially set to empty.
    \item[Step 5] Partition the elements of $H_{N,a}$ according to squares of
        side length $a^{1/2}$ by computing, for each $(x,y) \in H_{N,a}$,
        the values $m= \lfloor x/a^{1/2} \rfloor$ and
        $n= \lfloor y/a^{1/2} \rfloor$, and appending the point $(x,y)$ to Bin[m,n].
    \item[Step 6] Compute the $\phi(a-1)$ elements of $H_{N,a-1}$. For each $(x_1,y_1) \in H_{N,a-1}$:
        \begin{itemize}
            \item[Step 6a] Determine which bin it corresponds to by computing
            $m= \lfloor x_1/a^{1/2} \rfloor$ and $n= \lfloor y_1/a^{1/2} \rfloor$.
            \item[Step 6b] Loop through the nearby points $(x_0,y_0)$ of $H_{N,a}$
                whose coordinates lie, left and downwards, within $a^{1/2}$. Typically, this
                entails examining the four bins Bin[$m-\epsilon_1,n-\epsilon_2$],
                where $\epsilon_1,\epsilon_2 \in \{0,1\}$. However, slight care
                is needed when crossing over an edge of the $a\times a$ square- one should
                wrap to the opposite side of the square.
            \item[Step 6c] Set $\mu_0=x_0$, $\nu_0=y_0$, $\mu_1 = x_1-\mu_0$,
                $\nu_1 = y_1-\nu_0$, and check whether $(\mu_1 a+\mu_0)(\nu_1 a
                + \nu_0) = N$. If so, we have determined a non-trivial factor
                of $N$ and quit.
        \end{itemize}
\end{itemize}
\end{algorithm}

The storage requirement of $O(N^{1/3})$ can be improved to $O(N^{1/6})$ by generating the
solutions $(x,y)$ to $xy = N \mod a-\delta$ lying in one vertical strip of width $O(a^{1/2})$ at a time
(easy to do since we can choose $x$ as we please, which then determines
$y$). In general, we are then no longer free to generate all modular inverses at once, and must
compute inverses in intervals of size $a^{1/2}$, one at a time, at a cost, using the Euclidean 
algorithm, of $O(a^\epsilon)$ per inverse.

\subsection{Variant. $1<U\leq V<N$ without restriction.}
\label{sec:variant}

Say $U=N^\alpha$, $V=N^{1-\alpha}$, with $1/3 < \alpha \leq 1/2$.
We may assume that $\alpha >1/3$, for, if not, we can find $U$
by performing $O(N^{1/3})$ trial divisions.

Let $a=\lceil 2 N^{1/3} \rceil$ (we do, here, mean $2N^{1/3}$, rather than
$(2N)^{1/3}$ of the previous section, as explained below).

Instead of working with small squares of side length $a^{1/2}$, partition the
$a \times a$ square into rectangles of width $w$ and height $h$, with
$wh=N^{1/3}$. We would like to select $w$ roughly equal to $N^{\alpha-1/3}$ and
hence $h=N^{1/3}/w$ roughly equal to $N^{2/3-\alpha}$.

These rough values of $w$ and $h$ are needed to make sure that, using the same
notation as before, $(u_0,v_0)$ and $(u_0+u_1,v_0+v_1)$ are in the same, or
neighbouring, rectangles. More precisely, say $N^{\alpha-1/3} < w \leq
2N^{\alpha-1/3}$. Then $h=N^{1/3}/w \geq N^{2/3-\alpha}/2$. Then,
in~\eqref{eq:linear}, $u_1 = \lfloor U/a \rfloor \leq N^\alpha/\lceil 2 N^{1/3}
\rceil \leq N^{\alpha-1/3}/2 <w$, and $v_1 = \lfloor V/a \rfloor \leq
N^{1-\alpha}/\lceil 2 N^{1/3} \rceil \leq N^{2/3-\alpha}/2 \leq h$. Thus the
$x$-coordinates of $(u_0,v_0)$ and $(u_0+u_1,v_0+v_1)$ are $<w$ apart and the
$y$-coordinates are $\leq h$ apart.

Since we do not, a priori know $\alpha$, we cannot simply set $w$ and $h$.
Instead, we use an exponentially increasing set of $w$'s, for example starting
with $w=2$, and, repeatedly applying the above procedure, each time
doubling the size of $w$, until $w > N^{\alpha-1/3}$ and one successfully factors $N$.

The area of each rectangle is $N^{1/3}$,
and of the $a \times a$ square is approximately $N^{2/3}$, so there are  $O(N^{1/3})$ rectangles
(at the top and right edges these will typically be truncated), and, on average,
each contains $O(1)$ solutions to $xy = N \mod a-\delta$. Running through each rectangle and its
immediate neighbours, checking all pairs of points in these rectangles suggests
$O(N^{1/3})$ operations are needed for a particular choice of $w$ and $h$.
Since we might have to repeat this a few times, doubling the size of $w$,
the overall running time gets multiplied by $O(\log N)$ which is $O(N^\epsilon)$.

In Section~\ref{sec:moment}, a running time equal to $O(N^{1/3+\epsilon})$ is proven.

The steps described in this section are summarized below.

\begin{algorithm}
    Let $N=UV$ be a positive integer, and assume that $1<U \leq V < N$, with $U,V \in \Z$.
\begin{itemize}
    \item[Step 1] Carry out trial division on $N$ up to $N^{1/3}$. If a non-trivial factor of $N$ is
    found quit.
    \item[Step 2] Set $a=\lceil 2N^{1/3} \rceil$.
    \item[Step 3]Use the Euclidean algorithm to compute $\gcd(N,a-\delta)$ for $\delta=0,1$.
        If either gcd is $>1$ then we have determined a non-trivial factor of $N$ and quit.
    \item[Step 4] Compute and store, in two arrays, all the points of $H_{N,a}$ and $H_{N,a-1}$.
    \item[Step 5] Set $j=0$. While we have not succeeded in finding a non-trivial factor of $N$:
    \item
       \begin{itemize}
            \item[Step 5a] Increment $j$ by 1 and set $w=2^j$ and $h=N^{1/3}/w$.
            \item[Step 5b] For $0 \leq m < a/w$ and  $0 \leq n < a/h$,
            initialize a doubly indexed array, `Bin', whose elements,
            Bin[m,n], will contain lists of points and be used to partition
            $H_{N,a}$. Each bin is initially set to empty.
            \item[Step 5c] Partition the elements of $H_{N,a}$ according to rectangles of
            width $w$ and height $h$ by computing, for each $(x,y) \in H_{N,a}$,
            the values $m= \lfloor x/w \rfloor$ and
            $n= \lfloor y/h \rfloor$, and appending the point $(x,y)$ to Bin[m,n].
            \item[Step 5d] For each $(x_1,y_1) \in H_{N,a-1}$:
                \begin{itemize}
                    \item[Step 5d1] Determine which bin it corresponds to by computing
                    $m= \lfloor x_1/w \rfloor$ and $n= \lfloor y_1/h \rfloor$.
                    \item[Step 5d2] Loop through the nearby points $(x_0,y_0)$
                        of $H_{N,a}$ whose coordinates lie, left and downwards,
                        within $w$ and $h$ respectively. Typically, this
                        entails examining the four bins
                        Bin[$m-\epsilon_1,n-\epsilon_2$], where
                        $\epsilon_1,\epsilon_2 \in \{0,1\}$. However, slight
                        care is needed when crossing over an edge of the
                        $a\times a$ square- one should wrap to the opposite
                        side of the square.
                    \item[Step 5d3] Set $\mu_0=x_0$, $\nu_0=y_0$, $\mu_1 = x_1-\mu_0$,
                        $\nu_1 = y_1-\nu_0$, and check whether $(\mu_1 a+\mu_0)(\nu_1 a
                        + \nu_0) = N$. If so, we have determined a non-trivial factor
                        of $N$ and quit.
                \end{itemize}
            \item[Step 5e] Free up the memory used by `Bin'.
       \end{itemize}
\end{itemize}
\end{algorithm}

\section{Towards a subexponential bound}
\label{sec:subexp}

The above algorithm exploits the fact that when $a$ is large, and $\delta$ is small,
the points with coordinates $(U,V) \mod a-\delta$ are close to one another.
In fact they lie equally spaced on a line with common horizontal difference
$u_1$, and vertical difference $v_1$.

An obvious thing to try is to reduce the size of $a$. However, as $a$ decreases,
$u_1$ and $v_1$ increase so that not only do the points $(u_0,v_0)$ and $(u_0+u_1,v_0+v_1)$
move far apart, the latter point soon falls far outside the square of side length $a$.

To fix this, one can view (\ref{eq:linear}) as the base $a$ expansion of $U$ and $V$.
When $a$ is smaller, one could instead use a polynomial expansion
\begin{eqnarray}
    \label{eq: polynomial}
    U= u_{d_1} a^{d_1} + \ldots + u_1 a + u_0, \quad  0 \leq u_i < a \notag \\
    V= v_{d_2} a^{d_2} + \ldots + v_1 a + v_0, \quad  0 \leq v_i < a,
\end{eqnarray}
with $u_{d_1} \neq 0$ and $v_{d_1} \neq 0$.
For simplicity in what follows, assume that the degrees of both polynomials are equal,
$d_1 = d_2 = d$, so that both $U$ and $V$ satisfy $a^d \leq U,V < a^{d+1}$.

A polynomial of degree $d$ is determined uniquely by $d+1$ values.
Imitating the approach in Section~\ref{sec:intro}, we evaluate $N \mod a-\delta$ 
for $d+1$ values of $\delta$. A natural choice might be $\delta=0, \pm 1, \pm 2, \ldots$, 
but, to keep our polynomial values positive, we consider non-negative values of $\delta$,
and, for good measure, take extra values, $\delta=0,1,2,\ldots 2d$ (by extra, we mean $\delta \leq 2d$ rather
than $d\ leq d$).
Now,
\begin{equation}
    N = UV = (u_d \delta^d + \ldots + u_1 \delta + u_0)(v_d \delta^d + \ldots + v_1 \delta + v_0) 
    \mod a-\delta.
\end{equation}
Since $0 \leq u_j < a$, we have
\begin{equation}
    u_d \delta^d + \ldots + u_1 \delta + u_0 < a \lambda(d,\delta)
\end{equation}
where
\begin{equation}
   \lambda(d,\delta) = \delta^d + \delta^{d-1} + \ldots + 1 = (\delta^{d+1}-1)/(\delta-1)
    \sim \delta^d, \quad \text{as $\delta \to \infty$}.
\end{equation}
and similarly for the $v_j$'s.

For each $\delta$ one lists all solutions $(x,y)$ to
\begin{equation}
    xy = N \mod a-\delta
\end{equation}
\begin{equation}
    0 < x,y < a \lambda(d,\delta).
\end{equation}
The number of points $(x,y)$ for a given $\delta$ is $\phi(a-\delta)$ per
$a \times a$ square, and hence, overall, equals
\begin{equation}
    \phi(a-\delta) \lambda(d,\delta)^2 = O(a (2d)^{2d}).
\end{equation}
We are again assuming that $\gcd(a-\delta,N)=1$, otherwise one easily pulls out
a factor of $N$.

We need a method to recognize the solutions that we seek
$(u_d \delta^d + \ldots + u_0, v_d \delta^d + \ldots  + v_0)$ hiding amongst 
all the $(x,y)$'s. This leads to the question:

Let $X>0$ and let $S_0,S_1,\ldots,S_{2d}$ be $2d+1$ sets of points $\in \mathbb{Z}^2$
all of whose coordinates are positive and $\leq X$. Assume that amongst these points
there exists $2d+1$ points, one from each $S_\delta$, whose coordinates are described
by polynomials $u(\delta),v(\delta) \in \mathbb{Z}[\delta]$ of degree $d$. More precisely,
for each
$0 \leq \delta \leq 2d$ there exists a point $(x_\delta, y_\delta) \in S_\delta$ such that
\begin{eqnarray}
     x_\delta = u(\delta) = u_d \delta^d + \ldots + u_0 \notag \\
     y_\delta = v(\delta) = v_d \delta^d + \ldots + v_0.
\end{eqnarray}
Can one find these $2d+1$ points much more efficiently than by exhaustively searching 
through all possible $2d+1$ tuples of points? For example, can one find these
points in time $O(X^\alpha d^{\beta d})$ for some $\alpha,\beta>0$?

In our application, $X= O(a (2d)^{2d})$.
Since $N=UV$ and $a^d \leq U <V < a^{d+1}$, we have
$a< N^{1/(2d)}$. Assuming that there is an $O(X^\alpha d^{\beta d})$ time
algorithm for finding points with polynomial coordinates,
on taking $d$ proportionate to
\begin{equation}
    \left( \frac{\log N}{\log \log N} \right)^{1/2}
\end{equation}
one gets a factoring algorithm requiring
\begin{equation}
    \exp\left( \gamma (\log N \log \log N)^{1/2} \right)
\end{equation}
time and storage, for some $\gamma >0$.

One can cut back a bit on the search space, by noting, for example,
that the coefficients of $u(\delta)$ and $v(\delta)$ are integers
(this imposes a divisibility restriction
on finite differences between points lying on the polynomial), and,
in our particular application, that the coefficients are non-negative and bounded,
and this restricts the rate of growth of the polynomials. However, to get down to
a running time polynomial in $X$, one needs to do much better.


\section{Uniform distribution}
\label{sec:dist}

Let $\gcd(a,N)=1$. A classic application of Kloosterman sums shows that the points
$(x,y) \mod a$ satisfying $xy = N \mod a$ become uniformly distributed in the square of
side length $a$ as $a \to \infty$. While the tools used in this section
are fairly standard, they will also be applied in the next section to estimate the
running time of the Hide and Seek algorithm. Similar theorems can be found in the literature
\cite{BK} \cite{BCZ} \cite{FK} \cite{H} \cite{S} \cite{We}, often with restrictions to prime values of 
$a$ or to $N=1$.

Consider the following identity which detects pairs of integers $(x,y)$ such that
$xy = N \mod a$:
\begin{equation}
    \frac{1}{a} \sum_{k=0}^{a-1} e\left( \frac{k}{a} (y-\bar{x} N) \right)
    =\begin{cases}
          1 \quad \text{if $xy = N \mod a$} \\
          0 \quad \text{otherwise}
     \end{cases}
\end{equation}
where $e(z)= \exp(2 \pi i z)$, and where $\bar{x}$ stands for any integer
congruent to $x^{-1} \mod a$, if the inverse exists. Recall that we have assumed $\gcd(a,N)=1$
so that any solution to $xy=N \mod a$ must have $\gcd(x,a)=1$. Thus, for such solutions,
$x^{-1} \mod a$ exists.

Let $R$ be the rectangle bounded horizontally by $x_1,x_2 \in \mathbb{Z}$ and
vertically by $y_1, y_2 \in \mathbb{Z}$, where $0 \leq x_1 < x_2 \leq a$
and $0 \leq y_1 < y_2 \leq a$:
\begin{equation}
    \label{eq:rectangle}
    R = R(x_1,x_2,y_1,y_2) =
    \{ (x,y) \in \mathbb{Z}^2 | x_1 \leq x < x_2, y_1 \leq y < y_2 \}.
\end{equation}

Let $c_R(N,a)$ denote the number of pairs of integers $(x,y)$ that lie in the
rectangle $R$, and satisfy $xy = N \mod a$:
\begin{equation}
     c_R(N,a) =
     \sum_{{(x,y) \in R}\atop {xy = N \mod a}} 1.
\end{equation}
The identity above gives
\begin{equation}
     \label{eq:c_R}
     c_R(N,a) = \frac{1}{a} \sum_{k=0}^{a-1}
     \sum_{{(x,y) \in R}\atop {\gcd(x,a)=1}}
     e\left( \frac{k}{a} (y-\bar{x} N) \right).
\end{equation}
Notice that we only need to restrict $x$ to $\gcd(x,a)=1$ and that
$y$ runs over all residues in $y_1 \leq y < y_2$. This will allow us to
deal with the sum over $y$ as a geometric series.

The $k=0$ term provides the main contribution
while the other terms can be estimated using bounds for Kloosterman sums.
We require two lemmas. The first considers the main contribution,
and the second bounds the remaining terms.

\begin{lemma}
    The $k=0$ term in (\ref{eq:c_R}) equals
    \begin{equation}
        \frac{\text{area$(R)$}}{a^2} \phi(a) + O(a^\epsilon)
    \end{equation}
    for any $\epsilon>0$.

    \begin{proof}
    The $k=0$ term is
    \begin{equation}
         \label{eq:main contribution c_R}
         \frac{1}{a}
         \sum_{{(x,y) \in R}\atop {\gcd(x,a)=1}} 1
         = \frac{(y_2-y_1)}{a} \sum_{{x_1 \leq x < x_2}\atop {\gcd(x,a)=1}} 1.
    \end{equation}
    Using the Mobius function we have
    \begin{eqnarray}
        \label{eq: mobius}
        \sum_{{x_1 \leq x < x_2}\atop {\gcd(x,a)=1}} 1 &=&
        \sum_{x_1 \leq x < x_2}
        \sum_{d|\gcd(x,a)} \mu(d)
        = \sum_{d|a} \mu(d)
        \sum_{x_1/d \leq x < x_2/d} 1
        \notag \\
        &=& \sum_{d|a} \mu(d) ((x_2-x_1)/d + O(1))
        = (x_2-x_1) \prod_{p | a} (1-1/p) + O(\tau(a)),
    \end{eqnarray}
    where $\tau(a)$ equals the number of divisors of $a$ and is $O(a^\epsilon)$ for
    any $\epsilon>0$.
    This implies that the $k=0$ contribution to $c_R(N,a)$ equals
    \begin{equation}
        \frac{\text{area$(R)$}}{a^2} \phi(a) + O((y_2-y_1)a^{-1+\epsilon})
    \end{equation}
    which gives the lemma.
    \end{proof}
\end{lemma}

The next lemma bounds the contribution of the
$k \geq 1$ terms in (\ref{eq:c_R}).
\begin{lemma}
    For any $\epsilon >0$ we have
    \begin{equation}
         \frac{1}{a} \sum_{k=1}^{a-1}
         \sum_{{(x,y) \in R}\atop {\gcd(x,a)=1}}
         e\left( \frac{k}{a} (y-\bar{x} N) \right)
         = O(a^{1/2+\epsilon}).
    \end{equation}
    \begin{proof}
    One can separate the sum over $y$ and evaluate it as a geometric series
    obtaining for the lhs above
    \begin{equation}
         \frac{1}{a} \sum_{k=1}^{a-1}
         \frac{e\left( \frac{k}{a} y_2\right) - e\left( \frac{k}{a} y_1\right)}
         {e\left( \frac{k}{a}\right) -1}
         \sum_{{x_1 \leq x < x2}\atop {\gcd(x,a)=1}}
         e\left( \frac{-k}{a} \bar{x} N \right).
    \end{equation}

    Taking absolute values we get an upper bound of
    \begin{equation}
         \label{eq:upper bound 1}
         \frac{1}{a} \sum_{k=1}^{a-1}
         \frac{\left| \sin\left( \frac{\pi k}{a} (y_2-y_1\right) \right|}
         {\left| \sin\left( \frac{\pi k}{a} \right ) \right|}
         \left|
         \sum_{{x_1 \leq x < x2}\atop {\gcd(x,a)=1}}
         e\left( \frac{-k}{a} \bar{x} N \right)\right|.
    \end{equation}
    Next, notice that the terms $k$ and $a-k$ give the same contribution,
    so we may restrict our attention to just the terms $1 \leq k \leq (a-1)/2$.
    If $a-1$ is odd, the middle term is left out at a cost of $O(1)$, and the
    bound becomes
    \begin{equation}
         \label{eq:upper bound 2}
         \frac{2}{a} \sum_{1\leq k \leq (a-1)/2}
         \frac{\left| \sin\left( \frac{\pi k}{a} (y_2-y_1\right) \right|}
         {\left| \sin\left( \frac{\pi k}{a} \right ) \right|}
         \left|
         \sum_{{x_1 \leq x < x2}\atop {\gcd(x,a)=1}}
         e\left( \frac{-k}{a} \bar{x} N \right)\right| + O(1).
    \end{equation}

    The second sum above over $x$ can be expressed in terms of Kloosterman sums,
    and using estimates for Kloosterman sums one has
    \begin{equation}
        \label{eq: inequality kloosterman}
        \sum_{{x_1 \leq x < x2}\atop {\gcd(x,a)=1}}
        e\left( \frac{-k}{a} \bar{x} N \right)
        = O(a^{1/2+\epsilon} \gcd(k,a)^{1/2}).
    \end{equation}
    For a proof, see Lemma 4 on page 36 of Hooley's book \cite{Ho} where a proof is given
    (his $r$ corresponds to our $a$, and his $l$ is $-kN$.
    Also recall that we are assuming $\gcd(N,a)= 1$ so that $N$ does not appear
    in the gcd of the $O$ term).

    Furthermore, using the Taylor expansion of $\sin(x)$ one obtains the
    two inequalities
    \begin{eqnarray}
         \label{eq:inequality sin}
         \sin(x) \leq \min(x,1), \quad  x \geq 0, \notag \\
         1/\sin(x) < 2/x, \quad  0 < x < \pi/2.
    \end{eqnarray}
    For the second inequality, use $ 0 < x/2 < x - x^3/3! < \sin(x)$ in the
    stated interval.

    Applying (\ref{eq:inequality sin}) and (\ref{eq: inequality kloosterman})
    gives an upper bound for (\ref{eq:upper bound 2}) of
    \begin{equation}
         \label{eq:upper bound 3}
         O \left(
             a^{-1/2+\epsilon} 
             \sum_{1\leq k \leq (a-1)/2}
             \min(\frac{\pi k}{a} (y_2-y_1),1) \frac{2a}{\pi k}
             \gcd(k,a)^{1/2}
             + 1
         \right).
    \end{equation}
    Breaking up the sum into $1 \leq k \leq a/(\pi(y_2-y_1))$ and 
    $a/(\pi(y_2-y_1)) < k \leq (a-1)/2$, the sum over $k$ in the $O$
    term equals
    \begin{equation}
        \label{eq:two sums}
        2(y_2-y_1) \sum_{1\leq k \leq a/(\pi(y_2-y_1))} \gcd(k,a)^{1/2} +
        \frac{2a}{\pi} \sum_{a/(\pi(y_2-y_1)) < k \leq (a-1)/2} \gcd(k,a)^{1/2}/k.
    \end{equation}
    Both kinds of sums can be easily handled (the first can also be found in 
    Hooley). Let $X>0$. Then,
    \begin{equation}
        \label{eq:bound sum1}
        \sum_{1 \leq k \leq X} \gcd(k,a)^{1/2} 
        \leq \sum_{d|a} d^{1/2} \sum_{{1 \leq k \leq X }\atop{d|k}} 1
        \leq X \sum_{d|a} d^{-1/2} = O(X a^{\epsilon}).
    \end{equation}
    Next, let $0 <X_1 <X_2$. Then
    \begin{eqnarray}
        \sum_{X_1 < k \leq X_2} \gcd(k,a)^{1/2}/k
        \leq \sum_{d|a} d^{1/2} \sum_{{X_1 < k \leq X_2}\atop{d|k}} 1/k \notag \\
        = O\left(\log(X_2-X_1+2) \sum_{d|a} d^{-1/2} \right)
    \end{eqnarray}
    which equals
    \begin{equation}
        \label{eq:bound sum2}
        O(\log(X_2-X_1+2) a^\epsilon).
    \end{equation}

    Applying (\ref{eq:bound sum1}) and (\ref{eq:bound sum2}) to (\ref{eq:two sums}), 
    we have that (\ref{eq:upper bound 3}) is
    \begin{equation}
         O(a^{1/2+\epsilon}),
    \end{equation}
    completing the proof.
    \end{proof}
\end{lemma}

These two lemmas together give the following theorem.
\begin{theorem}
    Let $\gcd(N,a)=1$ and $R$ as described in (\ref{eq:rectangle}).
    Then, $c_R(N,a)$,
    the number of solutions $(x,y)$ to $xy = N \mod a$ with $(x,y)$ lying in the rectangle $R$,
    is equal to
    \begin{equation}
        \frac{\text{area$(R)$}}{a^2} \phi(a)
        + O(a^{1/2+\epsilon})
    \end{equation}
    for any $\epsilon>0$.
\end{theorem}
This theorem shows that the points $(x,y)$ satisfying $xy = N \mod a$
are uniformly dense in the sense that the rectangle $R$ contains its
fair share of solutions, so long as the area of $R$ is of larger size
than $a^{3/2+\epsilon}$. 

For example, if $R$ is a square, it needs to have side length at least
$a^{3/4+\epsilon}$ to contain its fair share of points. This is considerably
larger than the side length of $a^{1/2}$ that is used in the algorithm of Section 1.

The paper of Shparlinski \cite{S} contains many references to the
problem of uniform distribution and discusses improved results on average over $N$.

\section{Second moment and running time}
\label{sec:moment}

We now examine the assertion made in
Section 1 that $O(a^{1+\epsilon})$ time is needed to scan across all $a$ squares of side length
$a^{1/2}$ and their immediate neighbours, comparing all pairs of points contained in
said squares.

The running time of the algorithm in Section~\ref{sec:variant}, i.e. in the
case $U\leq V<2U$ depends on how the solutions to $xy = N \mod a$ and $x'y'  = N
\mod (a-1)$ are distributed amongst the small squares of side length $a^{1/2}$.
In Section~\ref{sec:runtime variant} we will consider the running time of the variant in
Section~\ref{sec:variant} which is used for the general situation $1< U\leq V< N$.

Let $S$ denote one such square, i.e. of side length $a^{1/2}$.
Then the running time needed to examine just the square $S$, looking at all pairs of
points $(x,y)$, $(x',y')$ in $S$ is $O(c_S(N,a) c_S(N,a-1))$, which,
by the arithmetic geometric inequality is $O(c_S(N,a)^2+ c_S(N,a-1)^2)$. The algorithm
also requires us to compare points in neighbouring squares, say $S_1$ and $S_2$,
which, similarly, takes $O(c_{S_1}(N,a)^2+ c_{S_2}(N,a-1)^2)$ time. Hence, the overall
running time to compare pairs of points is
\begin{equation}
    O \left( \sum_S c_S(N,a)^2+ c_S(N,a-1)^2\right), 
\end{equation}
the sum being over the roughly $a$ squares of side length $a^{1/2}$ that partition the $a \times a$
square $\{ (x,y) \in \mathbb{Z}^2 | 0 \leq x,y < a \}$ (at the top and right edges we get rectangles,
unless $a^{1/2}$ is an integer).

Consider now the contribution from the points $\mod a$:
\begin{equation}
    \label{eq: 2nd moment}
    \sum_S c_S(N,a)^2.
\end{equation}
For convenience, rather than deal with squares $S$ of side length $a^{1/2}$, we will estimate
~\eqref{eq: 2nd moment} by making a small
adjustment and partitioning the $a \times a$ square into squares of side length
\begin{equation}
    \label{eq:b}
    b=\lceil a^{1/2} \rceil.
\end{equation}
We also assume that $\gcd(b,a)=1$. If not, replace $b$ with $b+1$ until this condition holds.
By equation~(\ref{eq: mobius}), this will not take long to occur, so that, for any $\epsilon>0$,
$b=a^{1/2}+O(a^\epsilon)$.

Thus, consider the squares
\begin{equation}
    \label{eq: B}
    B_{ij} = 
    \{ (x,y) \in \mathbb{Z}^2 | i b  \leq x < (i+1) b , j b  \leq y < (j+1)b \}
\end{equation}
with $0 \leq i,j < a/b-1$.

Since $b \nmid a$ , these will not entirely cover the $a \times a$ square, but the number of points
$(x,y) \in \mathbb{Z}^2$
satisfying $xy = N \mod a$ that are neglected at the right most and top portions of the $a \times a$
square is, by~(\ref{eq: mobius}), $O(\phi(a) b/a)$, and these therefore contribute
$O(\phi(a)^2 b^2/a^2) = O(\phi(a))$ to~(\ref{eq: 2nd moment}).

The points $(x,y) \in \mathbb{Z}^2$ belonging to an $a^{1/2} \times a^{1/2}$ square $S$
are contained entirely in at most four squares, say $B_{i_1j_1}, B_{i_2j_2}, B_{i_1j_3}, B_{i_4j_4}$, 
of side length $b$. 
Therefore, 
\begin{equation}
    c_S(N,a)^2 \leq (c_{B_{i_1j_1}}(N,a)+c_{B_{i_2j_2}}(N,a)+c_{B_{i_3j_3}}(N,a)+c_{B_{i_4j_4}}(N,a))^2 
\end{equation}
which, by the Cauchy Schwartz inequality is
\begin{equation}
   \leq 4 (c_{B_{i_1j_1}}(N,a)^2+c_{B_{i_2j_2}}(N,a)^2+c_{B_{i_3j_3}}(N,a)^2+c_{B_{i_4j_4}}(N,a)^2).
\end{equation}
Since each $B$ square overlaps with $O(1)$ $S$ squares, we thus have that
\begin{equation}
    \sum_S c_S(N,a)^2 = O\left(\phi(a)+\sum_B c_B(N,a)^2\right),
\end{equation}
the $\phi(a)$ accounting for the contribution from the neglected portion at the right most and 
top portions of the $a \times a$ square.

A similar consideration for the points satisfying $xy = N \mod (a-1)$, 
partitioning the larger $a \times a$ square into
squares $D$ of side length $d$, where $d$ is the smallest integer greater than
$\lceil (a-1)^{1/2} \rceil$ which is coprime to $a$, gives the same kind of sum
\begin{equation}
    \sum_S c_S(N,a-1)^2 = O\left(\phi(a-1)+\sum_D c_D(N,a-1)^2\right).
\end{equation}

Therefore, we need to estimate the second moment
\begin{equation}
    \label{eq:2nd moment B}
    \sum_B c_B(N,a)^2
\end{equation}
where $B$ ranges over all $\lfloor a/b \rfloor^2$ squares of the form~(\ref{eq: B}).
To prove that the running time of the hide and seek algorithm of Section~\ref{sec:intro} is $O(N^{1/3+\epsilon})$ 
we need to prove that~(\ref{eq:2nd moment B}) is $O(a^{1+\epsilon})$.

\begin{theorem}
    Let
    \begin{equation}
         P= \{B_{ij}\}_{0 \leq i,j < a/b-1}
    \end{equation}
    Then 
    \begin{equation}
        \label{eq:thm 2nd moment}
        \sum_{B \in P} c_B(N,a)^2 = O(a^{1+\epsilon}).
    \end{equation}
\end{theorem}
\begin{proof}
Rather than look at just the $a \times a$ square, it is helpful to consider the
$ba \times ba$ square $\{ (x,y) \in \mathbb{Z}^2 | 0 \leq x,y < ba\}$.
The advantage of looking at the larger square will become apparent when we turn to the discrete
Fourier transform, and will be summing over all the $a$th roots of unity.

This larger square can be
partitioned into $b^2$ squares of side length $a$.
Because the solutions to $xy = N \mod a$ repeat $\mod a$, we can count each $c_B(N,a)^2$
once per $a \times a$ square, by summing $c_{B'}(N,a)^2=c_B(N,a)^2$ over all $b^2$ translates
$B'=B+(r_1 a, r_2 a)$ of $B$, with $0 \leq r_1,r_2 < b$.

On the other hand, we can also partition the $ba \times ba$ square into
$a^2$ squares of side length $b$:
\begin{equation}
    P_2 = \{B_{ij}\}_{0 \leq i,j \leq a-1}
\end{equation}
with $B_{ij}$ given by (\ref{eq: B}).

Each translate of a $B$ square, $B'=B+(r_1 a, r_2 a)$, is covered by at most four
$B_{ij} \in P_2$, and each $B_{ij} \in P_2$ overlaps at most four such translates of $B$.

Hence, applying the Cauchy-Schwartz inequality as before,
\begin{equation}
    \label{eq:with translates}
    b^2 \sum_{B \in P} c_B(N,a)^2 = O\left( \sum_{B \in P_2} c_B(N,a)^2 \right).
\end{equation}

To study $c_B(N,a)^2$ we multiply equation~(\ref{eq:c_R}) by its conjugate, giving
\begin{equation}
    c_B(N,a)^2 =
    \frac{1}{a^2}
    \sum_{0 \leq k_1, k_2 \leq a-1}
    \sum_{{(x_1,y_1) \in B}\atop {\gcd(x_1,a)=1}}
    \sum_{{(x_2,y_2) \in B}\atop {\gcd(x_2,a)=1}}
    e\left( \frac{k_1}{a} (y_1-\bar{x_1} N) - \frac{k_2}{a} (y_2-\bar{x_2} N)\right).
\end{equation}
Next, sum over all $B_{ij} \in P_2$, and
break up each sum over $(x,y)\in B_{ij}$ into a double sum $ib \leq x < (i+1)b$,
$jb \leq y < (j+1)b$,
\begin{eqnarray}
    \sum_{B \in P_2} c_B(N,a)^2 =
    \frac{1}{a^2} 
    \sum_{0 \leq k_1, k_2 \leq a-1}
    &
    \left(
        \displaystyle\sum_{i=0}^{a-1}
        \displaystyle\sum_{{ib \leq x_1,x_2 < (i+1)b }\atop {\gcd(x_1,a)=\gcd(x_2,a)=1}}
        e\left(-\frac{N}{a} (k_1\bar{x_1} - k_2\bar{x_2})\right)
    \right) \notag \\
    \label{eq:multi sum}
    \times &
    \left(
        \displaystyle\sum_{j=0}^{a-1}
        \displaystyle\sum_{jb \leq y_1,y_2 < (j+1)b }
        e\left(\frac{k_1y_1 - k_2y_2}{a}\right)
    \right).
\end{eqnarray}
Now, the inner most sum,
\begin{equation}
    \displaystyle\sum_{jb \leq y_1,y_2 < (j+1)b }
     e\left(\frac{k_1y_1 - k_2y_2}{a}\right),
\end{equation}
is a product of two geometric series and equals
\begin{equation}
     \label{eq:y geom series}
     e((k_1 - k_2)jb/a)
    \frac{e(k_1 b/a)-1}{e(k_1/a)-1}
    \frac{e(-k_2 b/a)-1}{e(-k_2/a)-1}.
\end{equation}
We understand $\frac{e(k b/a)-1}{e(k/a)-1}$ to equal $b$ if $k=0 \mod a$.
Summing (\ref{eq:y geom series})  over $0 \leq j \leq a-1$ gives
\begin{equation}
     \begin{cases}
        a \left| \frac{e(k_1 b/a)-1}{e(k_1/a)-1} \right|^2 \quad \text{if $k_1 = k_2 \mod a$} \notag \\
        0 \quad \text{otherwise} 
     \end{cases}
\end{equation}
(recall that we have chosen $b$ so that $\gcd(b,a)=1$). Therefore, only the terms with $k_1=k_2$
contribute to~(\ref{eq:multi sum}) and it equals
\begin{equation}
    \label{eq: simplified}
    \frac{1}{a}
    \sum_{k=0}^{a-1} 
    \left(
        \displaystyle\sum_{i=0}^{a-1}
        \displaystyle\sum_{{ib \leq x_1,x_2 < (i+1)b }\atop {\gcd(x_1,a)=\gcd(x_2,a)=1}}
        e\left(-\frac{N}{a} (k(\bar{x_1} - \bar{x_2}))\right)
    \right)
    \left| \frac{e(k b/a)-1}{e(k/a)-1} \right|^2.
\end{equation}
The $k=0$ term gives, on separating the sum over $x_1$ and $x_2$,
\begin{equation}
    \label{eq: 0 contribution}
    \frac{b^2}{a}
        \sum_{i=0}^{a-1}
        \left(
            \sum_{{ib \leq x < (i+1)b }\atop {\gcd(x,a)=1}} 1
        \right)^2
\end{equation}
which, by (\ref{eq: mobius}) and using $b \sim a^{1/2}$, equals
\begin{equation}
    \label{eq: k=0 contribution double sum}
    b^2 (\phi(a)b/a +O(a^\epsilon))^2 = O(\phi(a)^2).
\end{equation}

Next, we deal with the terms $1 \leq k \leq a-1$. The sum over $i$ in (\ref{eq: simplified}) equals
\begin{equation}
        \label{eq: i sum}
        \displaystyle\sum_{i=0}^{a-1}
        \displaystyle\sum_{ib \leq x_1,x_2 < (i+1)b }
        A_{x_1,x_2}(-Nk),
\end{equation}
where
\begin{equation}
        A_{x_1,x_2}(t)=
        \begin{cases}
            0 \quad \text{if $\gcd(x_1 x_2,a)>1$} \notag \\
            e(t(\bar{x_1} - \bar{x_2})/a) \quad \text{otherwise.}
        \end{cases}
\end{equation}

To analyze this sum, we use the two dimensional discrete Fourier transform
\begin{equation}
        \hat{A}_{m_1,m_2}(t)=
        \displaystyle\sum_{0 \leq x_1,x_2 \leq a-1}
        A_{x_1,x_2}(t) e\left(-\frac{m_1 x_1 + m_2 x_2}{a}\right),
\end{equation}
so that
\begin{equation}
        A_{x_1,x_2}(t)=
        \frac{1}{a^2}
        \displaystyle\sum_{0 \leq m_1,m_2 \leq a-1}
        \hat{A}_{m_1,m_2}(t) e\left(\frac{m_1 x_1 + m_2 x_2}{a}\right),
\end{equation}
and (\ref{eq: i sum}) equals, on changing order of summation,
\begin{equation}
        \label{eq:with ft}
        \frac{1}{a^2}
        \displaystyle\sum_{0 \leq m_1,m_2 \leq a-1}
        \hat{A}_{m_1,m_2}(-Nk) 
        \left(
        \displaystyle\sum_{0 \leq i \leq a-1}
        \displaystyle\sum_{ib \leq x_1,x_2 < (i+1)b }
        e\left(\frac{m_1 x_1 + m_2 x_2}{a}\right)
        \right).
\end{equation}
The bracketed sum over $i$ is similar to the sum over $j$ worked out above and equals
\begin{equation}
     \begin{cases}
        a \left| \frac{e(m_1 b/a)-1}{e(m_1/a)-1} \right|^2 \quad \text{if $m_2=-m_1 \mod a$} \notag \\
        0 \quad \text{otherwise}.
     \end{cases}
\end{equation}
Therefore, (\ref{eq:with ft}) equals
\begin{equation}
    \frac{1}{a}
    \sum_{m=0}^{a-1}
    \hat{A}_{m,a-m}(-Nk)
    \left| \frac{e(m b/a)-1}{e(m/a)-1} \right|^2.
\end{equation}

So, (\ref{eq: simplified}), and hence (\ref{eq:multi sum}), equals
\begin{equation}
    \label{eq: more simplified}
    \frac{1}{a^2}
    \sum_{k=0}^{a-1}
    \sum_{m=0}^{a-1}
    \hat{A}_{m,a-m}(-Nk)
    \left| \frac{e(m b/a)-1}{e(m/a)-1} \right|^2
    \left| \frac{e(k b/a)-1}{e(k/a)-1} \right|^2.
\end{equation}

But,
\begin{equation}
    \hat{A}_{m,a-m}(-Nk)=
    \sum_{{0 \leq x_1,x_2 \leq a-1}\atop{\gcd(x_1 x_2,a)=1}}
     e\left(-\frac{Nk}{a}(\bar{x_1} - \bar{x_2})\right)
     e\left(-\frac{m x_1 - m x_2}{a}\right)
     =
     \left|
         \sum_{{0 \leq x \leq a-1}\atop{\gcd(x,a)=1}}
         e\left(-\frac{Nk\bar{x} + m x}{a}\right)
     \right|^2.
\end{equation}
However, the sum on the rhs is a Kloosterman sum
\begin{equation}
     \sum_{{0 \leq x \leq a-1}\atop{\gcd(x,a)=1}}
     e\left(-\frac{Nk\bar{x} + m x}{a}\right)
     = S(-m,-Nk,a)
\end{equation}
and are known \cite{W}~\cite{IK} to satisfy the bound
\begin{equation}
     \label{eq: kloosterman bound}
     | S(-m,-Nk,a) | \leq \tau(a) \gcd(m,k,a)^{1/2} a^{1/2} = O(a^{1/2+\epsilon}\gcd(k,a)^{1/2})
\end{equation}
(recall we are assuming that $\gcd(N,a)=1$ so that $N$ does not appear on the rhs of this inequality).
Applying this bound to $\hat{A}_{m,a-m}(-Nk)$, shows that~(\ref{eq: more simplified})
is
\begin{equation}
    \label{eq:almost there}
    O\left(
        \frac{a^{\epsilon}}{a}
        \sum_{k=1}^{a-1}
        \sum_{m=0}^{a-1}
        \gcd(k,a) 
        \left| \frac{e(m b/a)-1}{e(m/a)-1} \right|^2
        \left| \frac{e(k b/a)-1}{e(k/a)-1} \right|^2
        +\phi(a)^2
    \right).
\end{equation}
The $\phi(a)^2$ terms comes from the $k=0$ contribution, (\ref{eq: k=0 contribution double sum}).
We must isolate this term, otherwise the estimate below will be too large.

Separating sums gives
\begin{equation}
    \label{eq:almost there2}
    O\left(
        \frac{a^{\epsilon}}{a}
        \left(
            \sum_{k=1}^{a-1}
            \gcd(k,a) 
            \left| \frac{e(k b/a)-1}{e(k/a)-1} \right|^2
        \right)
        \left(
           \sum_{m=0}^{a-1}
           \left| \frac{e(m b/a)-1}{e(m/a)-1} \right|^2
        \right)
        +\phi(a)^2
    \right).
\end{equation}
Both sums can be bounded using the same approach as for~(\ref{eq:upper bound 1})
in the previous section, namely: combining terms $k$ and $a-k$ (similarly for the $m$ sum, but taking
the $m=0$ term alone),
breaking up the sum into the terms with $k\leq a/(\pi b) \sim a^{1/2}/\pi$ (respectively, $m$), applying
inequalities~(\ref{eq:inequality sin}), estimating the resulting sums, using
$b \sim a^{1/2}$, we find, for any $\epsilon > 0$, that~(\ref{eq:almost there2})
equals
\begin{equation}
    O(b^2 a^{1+\epsilon}).
\end{equation}
We have thus estimated the sum that appears on the rhs of~(\ref{eq:with translates}).
The sum that we wish to bound appears on the lhs of (\ref{eq:with translates}) but with an extra
factor of $b^2$. Hence, dividing the above by $b^2$ gives $O(a^{1+\epsilon})$ for the sum in
theorem.

\end{proof}
Remark: In certain cases, such as when $a=p^2$, with $p$ prime, one can improve the above estimate
for the second moment to $O(a)$ by taking $b=p$ and, for $x=jp+l$, with
$\gcd(l,p)=1$, using $\bar{x} = \bar{l}^2 (l-jp)$.

\subsection{Running time of the variant, for $1 < U \leq V <N$}
\label{sec:runtime variant}

Instead of partitioning the $a \times a$ square into smaller squares of side length $b \sim a^{1/2}$,
we partition it into rectangles $R$ of width $w<a$ and height $h<a$, where $w,h \in \mathbb{Z}$ and, for
convenience,
$\gcd(w,a)= \gcd(h,a) = 1$.

We partition the $a \times a$ square and also the larger $wa \times ha$ rectangle into smaller
rectangles $R$:
\begin{eqnarray}
    R &=& R_{ij} = \{ (x,y) \in \mathbb{Z}^2 | i w  \leq x < (i+1) w , j h  \leq y < (j+1)h \} \notag \\
    Q &=& \{R_{ij}\}_{{0 \leq i < a/w-1} \atop {0 \leq j < a/h-1}} \notag \\
    Q_2 &=& \{R_{ij}\}_{0 \leq i,j \leq a-1}. \notag
\end{eqnarray}

As in Section 4.1, we have
\begin{equation}
    \label{eq:wh}
    wh \sum_{R \in Q} c_R(N,a)^2 = O\left( \sum_{R \in Q_2} c_R(N,a)^2 \right).
\end{equation}
with $wh$ appearing on the lhs since the large $wa \times ha$ rectangle has that many copies of 
the $a \times a$ square.

Using the discrete Fourier transform, as before,
\begin{equation}
    \sum_{R \in Q_2} c_R(N,a)^2 =
    \frac{1}{a^2}
    \sum_{k=0}^{a-1}
    \sum_{m=0}^{a-1}
    |S(-m,-Nk,a)|^2
    \left| \frac{e(m w/a)-1}{e(m/a)-1} \right|^2
    \left| \frac{e(k h/a)-1}{e(k/a)-1} \right|^2.
\end{equation}
This useful identity expresses the second moment for the larger $wa \times ha$ rectangle
as a sum involving Kloosterman sums.

The $k=0$ term can be estimated as in~(\ref{eq: k=0 contribution double sum}) and asymptotically equals
\begin{equation}
   \label{eq:wh k=0 term}
   \frac{h^2 w^2}{a^2} \phi(a)^2.
\end{equation}

For the $k \geq 1 $ terms, we use bound~(\ref{eq: kloosterman bound}) to estimate the Kloosterman sums 
and separate the double sum above to get a contribution of
\begin{equation}
    O\left(
        \frac{a^{\epsilon}}{a}
        \left(
            \sum_{k=1}^{a-1}
            \gcd(k,a) 
            \left| \frac{e(k h/a)-1}{e(k/a)-1} \right|^2
        \right)
        \left(
           \sum_{m=0}^{a-1}
           \left| \frac{e(m w/a)-1}{e(m/a)-1} \right|^2
        \right)
    \right).
\end{equation}
The first sum is estimated to equal $O(\tau(a)ah)$ while the second sum is $O(aw)$,
giving, for $k \geq 1 $ a contribution of
\begin{equation}
    \label{eq:wh k geq 1}
    O(a^{1+\epsilon} wh)
\end{equation}
for any $\epsilon>0$. Putting (\ref{eq:wh k=0 term}) and (\ref{eq:wh k geq 1}) together,
then dividing the lhs of (\ref{eq:wh}) by $wh$ gives the following estimate for the second moment:

\begin{theorem}
Let $1 < w,h < a$, with $\gcd(w,a)=\gcd(h,a)=1$.
Then, using the notation above, we have an estimate for the second moment
that depends on the area $wh$ of the rectangles $R$:
\begin{equation}
    \label{eq:wh theorem}
    \sum_{R \in Q} c_R(N,a)^2 = 
    \begin{cases}
       O(a^{1+\epsilon}) \quad \text{if $wh = O(a^{1+\epsilon})$,}\\
       O(wh \phi(a)^2/a^2) \quad \text{if $wh \gg a^\lambda$ for some $\lambda > 1$.}\\
    \end{cases}
\end{equation}

\end{theorem}
Remark: if $\gcd(w,a)=\gcd(h,a)=1$ does not hold, one can bound the lhs of~\eqref{eq:wh theorem}
by comparing with the same kind of sum, but where $w$ and $h$ are incremented, as before,
by at most $O(a^\epsilon)$ until this gcd condition holds. So long as $w,h \gg a^\epsilon$ to begin with,
the estimates in the above theorem are unaffected.

In Section 1.2, our choice of $w$ and $h$ has $wh = O(a)$, and the estimate for the second 
moment is thus $O(a^{1+\epsilon})$, as in the previous section.

The second estimate of the theorem (not relevant for our particular application),
$O(wh \phi(a)^2/a^2)$, can probably be turned into an asymptotic
formula and a central limit theorem proven. This will remain
an inquiry for the future.

\subsubsection{Acknowledgements}

I wish to thank Andrew Granville, Carl Pomerance, and Matthew Young for helpful
feedback.

\bibliographystyle{amsplain}

\end{document}